\documentclass[11pt,draft]{amsart}
\usepackage{mathrsfs}
\usepackage{amssymb}
\usepackage{amsmath}
\usepackage{amsfonts}
\usepackage{amsfonts,enumerate}
\usepackage{pifont}
\usepackage{enumerate}
\usepackage{graphics}
\usepackage{verbatim}
\usepackage{amsthm}

\numberwithin{equation}{section}
\newtheorem{theorem}{Theorem}[section]
\newtheorem{corollary}{Corollary}[section]

\newtheorem{lemma}{Lemma}[section]

\newtheorem{remark}{Remark}[section]

\newcommand{\A}{{\mathcal A}}

\newcommand{\g}{\gamma}

\newcommand{\8}{\infty}
\newcommand{\el}{\ell}

\newcommand{\be}{\begin{eqnarray*}}
\newcommand{\ee}{\end{eqnarray*}}
\newcommand{\beq}{\begin{equation}}
\newcommand{\eeq}{\end{equation}}
\newcommand{\beqn}{\begin{equation*}}
\newcommand{\eeqn}{\end{equation*}}
\newcommand{\bs}{\begin{split}}
\newcommand{\es}{\end{split}}

\numberwithin{equation}{section}

\begin{document}
\title{Tent spaces and Littlewood-Paley $g$-functions associated with Bergman spaces\\ in the unit ball of $\mathbb{C}^n$}


\author{Zeqian Chen}

\address{Wuhan Institute of Physics and Mathematics, Chinese
Academy of Sciences, 30 West District, Xiao-Hong-Shan, Wuhan
430071,China}
\email{chenzeqian@hotmail.com}


\author{Wei Ouyang}

\address{Wuhan Institute of Physics and Mathematics, Chinese Academy of Sciences, 30 West District, Xiao-Hong-Shan, Wuhan 430071, China
and Graduate University of Chinese Academy of Sciences, Beijing 100049, China}


\date{}
\maketitle

\markboth{Z. Chen and W. Ouyang}%
{Bergman spaces}

\begin{abstract}
In this paper, a family of holomorphic spaces of tent type in the unit ball of $\mathbb{C}^n$ is introduced, which is closely related to maximal and area integral functions in terms of the Bergman metric. It is shown that these spaces coincide with Bergman spaces. Furthermore, Littlewood-Paley type $g$-functions for the Bergman spaces are introduced in terms of the radial derivative, the complex gradient, and the invariant gradient. The corresponding characterizations for Bergman spaces are obtained as well. As an application, we obtain new maximal and area integral characterizations for Hardy-Sobolev spaces.
\end{abstract}



\section{Introduction}\label{intro}

There is a mature and powerful real-variable theory for Hardy spaces of several complex variables which has distilled some of the essential oscillation and cancellation
behavior of holomorphic functions and then found that behavior ubiquitous. A good introduction to that is \cite{Rudin1980, Stein1972}; a more recent and fuller account is in \cite{AB1988, D1994, GL1978, GP2002, KL1995, KL1997} and references therein. However, the real-variable theory of the Bergman space is less well developed, even in the case of the unit disc (cf. \cite{DS2004}).

Recently, in \cite{CO} the present authors established maximal and area integral characterizations of Bergman spaces in the unit ball of $\mathbb{C}^n.$ The characterizations are in terms of maximal functions and area functions on Bergman balls involving the radial derivative, the complex gradient, and the invariant gradient. Moreover, the characterizations extend to cover Besov-Sobolev spaces. A special case of this is a characterization of $\mathrm{H}^p$ spaces involving only area functions on Bergman balls.

In this paper, we continue this line of investigation. We will consider the holomorphic spaces of tent type in the unit ball of $\mathbb{C}^n,$ which is closely related to maximal and area integral functions in terms of the Bergman metric. We show that these spaces actually coincide with Bergman spaces and hence obtain real-variable characterizations of them. Consequently, they provide the natural setting for the study of maximal and area integral functions for Bergman spaces in the unit ball of $\mathbb{C}^n.$ This motivation aries from the tent spaces in $\mathbb{R}^n,$ that were introduced and developed by Coifman, Meyer and Stein in \cite{CMS1985}. Our proofs involve some sharp estimates of the Bergman kernel function and Bergman metric. However,  since the Bergman metric $\beta$ is non-doubling, we will utilize some techniques of non-homogeneous harmonic analysis developed in \cite{NTV1998}. Furthermore, we will introduce Littlewood-Paley type $g$-functions for the Bergman spaces in terms of the radial derivative, the complex gradient, and the invariant gradient, which are variants of the ones used for Hardy spaces. The corresponding characterization are presented as well.

This paper is organized as follows. In Section \ref{notation} we present some notations and collect a number of auxiliary (and mostly elementary) facts about the Bergman metric and kernel functions. We state our main results in Section \ref{mainresult}. Section \ref{prooftentspace} is devoted to prove one of the main results that the holomorphic spaces of tent type coincide with Bergman spaces. To this end, the $L^p$-boundedness of non-central Hardy-Littlewood maximal function operators defined in terms of Bergman balls in $\mathbb{B}_n$ will be proved by the arguments of non-homogeneous harmonic analysis developed in \cite{NTV1998}. In Section \ref{proofg-funct}, we give the proof of another main result, that is, Littlewood-Paley type $g$-function characterizations of the Bergman spaces. Finally, in Section \ref{Besov} we give an application of our main results to obtaining new maximal and area integral characterizations of Besov spaces, including Hardy-Sobolev spaces.

In what follows, $C$ always denotes a constant depending (possibly) on $n, q, p, \g$ or $\alpha$ but not on $f,$ which may be
different in different places. For two nonnegative (possibly infinite) quantities $X$ and $Y,$ by $X \lesssim Y$ we mean that there
exists a constant $C>0$ such that $ X \leq C Y$ and by $X \thickapprox Y$ that $X \lesssim Y$ and $Y \lesssim X.$ Any notation and terminology not otherwise explained, are as used in \cite{Zhu2005} for spaces of holomorphic functions in the unit ball of $\mathbb{C}^n.$

\section{Basic facts and notation}\label{notation}

Throughout the paper we fix a positive integer $n \ge 1$ and a parameter $\alpha > -1.$ We denote by $\mathbb{C}^n$ the Euclidean space of complex dimension $n$. For $z=(z_1,\cdots,z_n)$ and $w=(w_1,\cdots,w_n)$ in $\mathbb{C}^n,$ we write
\be
\langle z,w\rangle=z_1\overline{w}_1+\cdots+z_n\overline{w}_n,
\ee
where $\overline{w}_k$ is the complex conjugate of $w_k.$ We also write
\be
|z|=\sqrt{|z_1|^2+\cdots+|z_n|^2}.
\ee
The open unit ball in $\mathbb{C}^n$ is the set
\be
{\mathbb{B}_n} =\{\,z\in {\mathbb{C}^n}: |z|<1\}.
\ee
The boundary of $\mathbb{B}_n$ will be denoted by $\mathbb{S}_n$ and is called the unit sphere in $\mathbb{C}^n,$ i.e.,
\be
{\mathbb{S}_n}=\{\,z\in {\mathbb{C}^n}: |z|=1\}.
\ee

For $\alpha$ $\in$ $\mathbb{R}$, the weighted Lebesgue measure $dv_{\alpha}$ on $\mathbb{B}_n$ is defined by
\be
dv_{\alpha}(z)=c_{\alpha}(1-|z|^2)^{\alpha}dv(z)
\ee
where $c_{\alpha}=1$ for $\alpha\le-1$ and $c_{\alpha}=\Gamma(n+\alpha+1)/ [n!\Gamma(\alpha+1)]$ if $\alpha>-1$, which is a normalizing constant so that $dv_{\alpha}$
is a probability measure on $\mathbb{B}_n.$ In the case of $\alpha=-(n+1)$ we denote the resulting measure by
\be
d\tau(z)=\frac{dv}{(1-|z|^2)^{n+1}},
\ee
and call it the invariant measure on $\mathbb{B}^n,$ since $d\tau=d\tau\circ\varphi$ for any automorphism $\varphi$ of $\mathbb{B}^n.$

For $\alpha > -1$ and $p>0$ the (weighted) Bergman space $\mathcal{A}^p_{\alpha}$ consists of holomorphic functions $f$ in $\mathbb{B}_n$ with
\be
\|f\|_{p,\,\alpha}=\left ( \int_{\mathbb{B}_n}|f(z)|^pdv_{\alpha}(z) \right )^{1/p}<\infty.
\ee
Thus,
\be
\mathcal{A}^p_{\alpha} = \mathcal{H} (\mathbb{B}_n) \cap L^p (\mathbb{B}_n, d v_{\alpha}),
\ee
where $\mathcal{H} (\mathbb{B}_n)$ is the space of all holomorphic functions in $\mathbb{B}_n.$ When $\alpha =0$ we simply write $\A^p$ for $\A^p_0.$ These are the usual Bergman spaces. Note that for $1 \le p < \8,$ $\mathcal{A}^p_{\alpha}$ is a Banach space under the norm $\|\ \|_{p,\,\alpha}.$ If $0 < p <1,$ the space $\mathcal{A}^p_{\alpha}$ is a quasi-Banach space with $p$-norm $\| f \|^p_{p, \alpha}.$

Recall that $D(z,\gamma)$ denotes the Bergman metric ball at $z$
\be
D(z, \gamma) = \{w \in \mathbb{B}_n\;: \; \beta (z, w) < \g \}
\ee
with $\gamma >0,$ where $\beta$ is the Bergman metric on $\mathbb{B}_n.$ It is known that
\be
\beta (z, w) = \frac{1}{2} \log \frac{1 + | \varphi_z (w)|}{1 - | \varphi_z (w)|},\quad z, w \in \mathbb{B}_n,
\ee
whereafter $\varphi _z$ is the bijective holomorphic mapping in $\mathbb{B}_n,$ which satisfies $\varphi _z (0)=z$, $\varphi _z(z)=0$ and
$\varphi _z\circ\varphi _z = id.$ If $\mathbb{B}_n$ is equipped with the Bergman metric $\beta,$ then $\mathbb{B}_n$ is a separable metric space. We shall call
$\mathbb{B}_n$ a separable metric space instead of $(\mathbb{B}_n,\beta).$

For the sake of convenience, we collect some elementary facts on the Bergman metric and holomorphic functions in the unit ball of $\mathbb{C}^n$ as follows.

\begin{lemma}\label{le:JfunctEsti} (cf. Theorem 1.12 in \cite{Zhu2005})
Let $c>0.$ Suppose $\alpha >-1.$ Then
\be
J_{c, \alpha}(z): = \int_{{\mathbb{B}}_n}\frac{d v_{\alpha}(w)}{|1-\langle z,w\rangle|^{n+1+ \alpha +c}} \thickapprox (1-|z|^2)^{-c}
\ee
as $|z|\rightarrow 1^-.$
\end{lemma}

\begin{lemma}\label{le:BergmanballVolumeEsti} (cf. Lemma 1.24 in \cite{Zhu2005})
For any real $\alpha$ and positive $\g$ there exists a constant $C_{\gamma}>0$ such that
\be
C_{\gamma}^{-1}(1-|z|^2)^{n+1+\alpha}\le v_{\alpha}(D(z,\g))\le C_{\gamma}(1-|z|^2)^{n+1+\alpha}
\ee
for all $z\in\mathbb{B}_n$.
\end{lemma}

\begin{lemma}\label{le:BergmanOperEsti} (cf. Lemma 2.10 in \cite{Zhu2005})
Let $a$ and $b$ be two real numbers. We define an integral operator $S$ by
\be
Sf(z)=(1-|z|^2)^a\int_{{\mathbb{B}}_n}\frac{(1-|w|^2)^b}{|1-\langle z,w\rangle|^{n+1+a+b}}f(w)dv(w).
\ee
Then for $-\infty<t<\infty$ and $1\leq p<\infty,$ $S$ is bounded on $L^p({\mathbb{B}}_n, dv_t)$ if and only if
$-pa<t+1<p(b+1).$
\end{lemma}

\begin{lemma}\label{le:BasicEsti01} (cf. Lemma 2.20 in \cite{Zhu2005})
For each $\gamma>0,$
\be
1-|a|^2 \approx 1-|z|^2 \approx |1-\langle a,z\rangle|
\ee
for all $a$ and $z$ in $\mathbb{B}_n$ with $\beta(a,z)<{\gamma}.$
\end{lemma}

\begin{lemma}\label{le:BergmanFunctEsti01} (cf. Lemma 2.24 in \cite{Zhu2005})
Suppose $\g >0, p>0,$ and $\alpha > -1.$ Then there exists a constant $C>0$ such that for any $f \in \mathcal{H} (\mathbb{B}_n),$
\be
|f(z) |^p \le \frac{C}{v_{\alpha}(D(z,\g))} \int_{D(z,\g)} | f(w)|^p d v_{\alpha} (w), \quad \forall z \in \mathbb{B}_n.
\ee
\end{lemma}

\begin{lemma}\label{le:BasicEsti02} (cf. (2.20) after the proof of Lemma 2.27 in \cite{Zhu2005})
For each $\gamma>0,$
\be
|1-\langle z,u\rangle| \approx |1-\langle z,v\rangle|
\ee
for all $z$ in $\bar{\mathbb{B}}_n$ and $u,v$ in $\mathbb{B}_n$ with $\beta(u,v)<\gamma.$
\end{lemma}

\begin{lemma}\label{le:BergmanFunctEsti02} (cf. Lemma 3.3 in \cite{Zhu2005})
Suppose $\beta$ is a real constant and $g\in L^1({\mathbb{B}}_n,dv)$. If
\be
f(z)=\int_{{\mathbb{B}}_n}\frac{g(w)dv(w)}{(1-\langle z,w\rangle)^{\beta}},\  z\in{\mathbb{B}}_n,
\ee
then
\be
\left|\tilde{\nabla}  f (z) \right|\leq\sqrt{2}|\beta|(1-|z|^2)^{\frac{1}{2}}\int_{{\mathbb{B}}_n}\frac{|g(w)|dv(w)}{|1-\langle z,w\rangle|^{\beta+\frac{1}{2}}}
\ee
for all $z\in{\mathbb{B}}_n$.
\end{lemma}

\begin{lemma}\label{le:Bergmanatomicdecomp}
Suppose $p>0, \alpha>-1,$ and $b> n \max \{1, 1/p \} + (\alpha+1)/p.$ Then there exists a sequence $\{a_k\}$ in $\mathbb{B}_n$ such that
$\mathcal{A}^p_{\alpha}$ consists exactly of functions of the form
\beq\label{eq:atomdecomp}
f(z)=\sum^{\infty}_{k=1}c_{k}\frac{(1-|a_k|^2)^{(pb-n-1-\alpha)/p}}{(1-\langle z,a_{k}\rangle)^b},\quad  z\in\mathbb{B}_n,
\eeq
where $\{c_k\}$ belongs to the sequence space $\el^p$ and the series converges in the norm topology of $\mathcal{A}^p_{\alpha}.$ Moreover,
\be
\int_{\mathbb{B}_n}|f(z)|^pdv_{\alpha}(z) \approx  \inf \Big \{ \sum_{k}|c_{k}|^p \Big \},
\ee
where the infimum runs over all decompositions of $f$ described above.
\end{lemma}

Lemma \ref{le:Bergmanatomicdecomp} is the atom decomposition for Bergman spaces due to Coifman and Rochberg \cite{CR1980} (see also \cite{Zhu2005}, Theorem 2.30).

\section{Statement of main results}\label{mainresult}

The basic functional used below is the one mapping functions in $\mathbb{B}_n$ to functions in $\mathbb{B}_n,$ given by
\beq\label{eq:TentFunct1}
A^{(q)}_{\gamma}(f)(z) = \left(\int_{D(z, \g)} | f (w)|^qd\tau(w)\right)^{\frac{1}{q}}
\eeq
if $1 < q < \infty,$ and
\beq\label{eq:TentFunct2}
A^{(\infty)}_{\gamma} (f) (z) = \sup_{w \in D(z, \g)} | f (w)|,\quad \mathrm{when} \ q=\infty.
\eeq
Then, given $0<p<\infty, 1 < q \le \8, $ and $\alpha>-1,\; \gamma > 0,$ the ``holomorphic space of tent type" $T^p_{q, \g} (\mathbb{B}_n, d v_{\alpha})$ in $\mathbb{B}_n$ is defined as the set of all $f \in \mathcal{H} (\mathbb{B}_n)$ so that $A^{(q)}_{\gamma}(f) \in L^p_{\alpha},$ equipped with a norm (or, quasi-norm)
\be
\|f \|_{T^p_{q, \g} (\mathbb{B}_n, d v_{\alpha})} = \|A^{(q)}_{\g}(f)\|_{p,\alpha}.
\ee
We will show that the space $T^p_{q, \g} (\mathbb{B}_n, d v_{\alpha})$ is independent of the choice of $\g,$ that is, $T^p_{q, \g_1} (\mathbb{B}_n, d v_{\alpha}) = T^p_{q, \g_2} (\mathbb{B}_n, d v_{\alpha})$ for any $0< \g_1, \g_2 < \8.$ We simply write $T^p_{q, \alpha} = T^p_{q, \g} (\mathbb{B}_n, d v_{\alpha})$ as follows. When $p\ge 1,$ $T^p_{q, \alpha}$ are Banach spaces under the norm $\|f \|_{T^p_{q, \alpha}}.$ If $0 < p < 1,$ the space $T^p_{q, \alpha}$ is a quasi-Banach space with $p$-norm $\|f \|_{T^p_{q, \alpha}}.$

The case of $q=\infty$ and $0<p<\infty$ was studied in \cite{CO}. Actually, the resulting tent type spaces $T^p_{\8, \alpha}$ are proved to be Bergman spaces $\mathcal{A}^p_{\alpha}.$

It is well known that the Hardy-Littlewood maximal function operator has played important role in harmonic analysis. To cater our estimates, we use a variant of the non-central Hardy-Littlewood maximal function operator acting on the weighted Lebesgue spaces $L^p_\alpha({\mathbb{B}_n})$, namely,
\beq\label{eq:MaxFunct2}
M^{(q)}_{\g} (f) (z) = \sup_{z \in D(w, \g)}\left(\frac{1}{v_{\alpha}(D(w,\g))}\int_{D(w,\g)}|f|^qdv_{\alpha}\right)^{\frac{1}{q}}
\eeq
for $1 \le q < \8.$ We simply write $M_{\g} (f) (z): = M^{(1)}_{\g} (f) (z).$

Now we are ready to formulate one of the main results of the present work.

\begin{theorem}\label{th:tentspace}
Suppose $\gamma >0, 1 \le q < \8,$ and $\alpha > -1.$ Let $0< p <\8.$ Then for any $f \in \mathcal{H} (\mathbb{B}_n),$ the following conditions are equivalent:
\begin{enumerate}[{\rm (1)}]

\item $f \in \A^p_{\alpha}.$

\item $A^{(q)}_{\g} (f)$ is in $L^p(\mathbb{B}_n, d v_{\alpha}).$

\item $M^{(q)}_{\g} (f)$ is in $L^p(\mathbb{B}_n, d v_{\alpha}).$

\end{enumerate}
Moreover,
\be
\| f \|_{\A^p_{\alpha}} \approx \| f \|_{T^p_{q, \alpha}} \approx \| M^{(q)}_{\g} (f) \|_{p, \alpha},
\ee
where the comparable constants depend only on $\gamma, q, \alpha, p,$ and $n.$
\end{theorem}

Note that the Bergman metric $\beta$ is non-doubling on $\mathbb{B}^n$ and so $(\mathbb{B}_n, \beta, d v_{\alpha})$ is a non-homogeneous space. The proof of the above theorem will involve some techniques of non-homogeneous harmonic analysis developed in \cite{NTV1998} (see Section \ref{prooftentspace} below).

In order to state the characterizations of Bergman spaces in terms of Littlewood-Paley type $g$-functions, we require some more notation. For any $f \in \mathcal{H} ( \mathbb{B}_n)$ and $z = (z_1, \ldots, z_n) \in \mathbb{B}_n$ we define
\be
\mathcal{R} f (z) = \sum^n_{k=1} z_k \frac{\partial f (z)}{\partial z_k}
\ee
and call it the radial derivative of $f$ at $z.$ The complex and invariant gradients of $f$ at $z$ are respectively defined as
\be
\nabla f(z) = \Big ( \frac{\partial f (z) }{\partial z_1},\ldots, \frac{\partial f (z) }{\partial z_n} \Big )\; \text{and}\; \widetilde{\nabla} f (z) = \nabla(f \circ \varphi_z)(0).
\ee

Given $1 < q < \8.$ We define for each $f \in \mathcal{H} (\mathbb{B}_n)$ and $z \in \mathbb{B}_n:$
\begin{enumerate}[{\rm (i)}]

\item The {\it radial} Littlewood-Paley $g$-function
\be
G^{(q)}_{\mathcal{R}}( f) (z) = \left ( \int_0^1   |(1-r) \mathcal{R} f (r z) |^q \frac{d r}{1-r} \right )^{\frac{1}{q}}.
\ee

\item The {\it complex gradient} Littlewood-Paley $g$-function
\be
G^{(q)}_{\nabla} ( f) (z) = \left ( \int_0^1  |(1-r) \nabla f (r z) |^q \frac{d r}{1-r} \right )^{\frac{1}{q}}.
\ee

\item The {\it invariant gradient} Littlewood-Paley $g$-function
\be
G^{(q)}_{\tilde{\nabla}} ( f) (z) = \left (  \int_0^1 |\tilde{\nabla} f (r z) |^q \frac{d r}{1- r|z|} \right )^{\frac{1}{q}}.
\ee

\end{enumerate}

We state another main result of this paper as follows.

\begin{theorem}\label{th:g-FunctCharat}
Suppose $1 < q < \8$ and $\alpha > -1.$ Let $0 < p < \8.$ Then, for any $f \in \mathcal{H} (\mathbb{B}_n)$ the following conditions are equivalent:
\begin{enumerate}[{\rm (1)}]

\item $f \in \mathcal{A}^p_{\alpha}.$

\item $G^{(q)}_{\mathcal{R}}(f)$ is in $L^p (\mathbb{B}_n, d v_{\alpha}).$

\item $G^{(q)}_{\nabla}( f)$ is in $L^p (\mathbb{B}_n, d v_{\alpha}).$

\item $G^{(q)}_{\tilde{\nabla}}( f)$ is in $L^p (\mathbb{B}_n, d v_{\alpha}).$

\end{enumerate}
Moreover, the following three quantities
\be
\| G^{(q)}_{\mathcal{R}}(f) \|_{p, \alpha},\; \| G^{(q)}_{\nabla}( f) \|_{p, \alpha}, \; \| G^{(q)}_{\tilde{\nabla}}( f) \|_{p, \alpha},
\ee
are all comparable to $\| f - f(0) \|_{p, \alpha},$ where the comparable constants depend only on $q, \alpha, p,$ and $n.$
\end{theorem}

\begin{remark}\label{rk:g-function}
Evidently, the results still hold true when we replace $1-r$ by $1-r |z|$ in the definitions of both $G^{(q)}_{\mathcal{R}}(f)$ and $G^{(q)}_{\nabla}( f).$ This is so because
\be
(1- |z|^2) | \mathcal{R} f (z) | \le (1- |z|^2) |\nabla f (z) | \le |\tilde{\nabla} f (z)|, \quad \forall z \in \mathbb{B}_n.
\ee
However, this is not the case if one uses $1-r$ instead of $1- r|z|$ in the definition of $G^{(q)}_{\tilde{\nabla}}( f).$ In fact,
\be
\int_0^1 |\tilde{\nabla} f (rz) |^2 \frac{d r}{1-r} = \8,\quad \forall z \in \mathbb{B}_n,
\ee
if $f(z)=z_1,$ but $f \in \mathcal{A}^p_{\alpha}$ for all $0 < p < \8$ and $\alpha >-1.$
\end{remark}

\section{Proof of Theorem \ref{th:tentspace}}\label{prooftentspace}

The goal of this section is to prove Theorem \ref{th:tentspace}. To this end, we need to prove the boundedness of non-central Hardy-Littlewood maximal function operator $M_{\g}$ as follows.

\begin{lemma}\label{le:Maximalfunction}
Let $\alpha > -1$ and $\g >0.$ The non-central Hardy-Littlewood maximal function operator $M_{\g}$ is bounded on $L^{p}_{\alpha}(\mathbb{B}_n)$ for each $1 < p \le \infty$ and acts from $
L^{1}_{\alpha}(\mathbb{B}_n)$ to $L^{1,\infty}_{\alpha},$ where
\be
L^{1, \8}_{\alpha} = \left \{ f:\; \| f \|_{L^{1, \8}_{\alpha}} \triangleq \sup_{\lambda > 0} \lambda v_{\alpha} \big ( \{ z \in \mathbb{B}_n :\; |M_{\g}f (z ) | > \lambda \} \big ) < \8 \right \}.
\ee
\end{lemma}

This can be obtained by using the following celebrated Vitali covering lemma, which can be founded in \cite{NTV1998}.

\begin{lemma}\label{le:Vitali}
Fix some $R>0$. Let $X$ be a separable metric space, $E$ any subset of $X,$ and let $\{B(x,r_x)\}_{x\in E}$ be a family of balls of radii $0<r_x<R$. Then there exists a countable subfamily
$\{B(x_j,r_j)\}_{j=1}^{\infty}$(where $x_j\in E$ and $r_{j} : =
 r_{x_j}$) of disjoint balls such that $E\subset\cup_j\{B(x_j,3r_j)\}_{j=1}^{\infty}.$
\end{lemma}

{\it Proof of Lemma \ref{le:Maximalfunction}}.\; The boundedness of $M_{\g}$ on $L^{\infty}_{\alpha}(\mathbb{B}_n)$ is obvious. Let $E_{\lambda}: = \left \{z \in \mathbb{B}_n:\; M_{\g}(f)(z)> \lambda \right\}.$ For each $z\in E_{\lambda},$ there exist some $w_z \in \mathbb{B}_n$ such that
\be
\int_{D(w_z, \g)}|f|dv_{\alpha}>\lambda v_{\alpha}(D(w_z, \g)).
\ee
Note that $z \in D(w_z, \g) \subset E_{\lambda}$ and $E_{\lambda} = \cup_{z\in E_{\lambda}} D(w_z, \g).$ Applying Lemma \ref{le:Vitali} and Lemma \ref{le:BergmanballVolumeEsti}, we choose the corresponding collection of pairwise disjoint Bergman metric ball $D(w_i, \g)$ and have
\be\begin{split}
v_{\alpha}(E_{\lambda}) & \le \sum_i v_{\alpha}(D(w_i,3 \g))\\
& \le C_{\g}\sum_{i} v_{\alpha}(D(w_i, \g)) \le \frac{C}{\lambda}\sum_i \int_{D(w_i, \g)} |f| d v_{\alpha}.
\end{split}\ee
This yields that
\be
v_{\alpha}(E_{\lambda})\le\frac{C}{\lambda}\int_{\mathbb{B}_n}|f|dv_{\alpha},
\ee
that is $M_{\g}$ is bounded from $L^{1}_{\alpha}(\mathbb{B}_n)$ to $L^{1,\infty}_{\alpha}.$

The boundedness of $M_{\g}$ on $L^{p}_{\alpha}(\mathbb{B}_n)$ for $1<p<\infty$ then follows from the Marcinkiewicz interpolation theorem.
\hfill$\Box$

\begin{remark}\label{rk:Maximalfunction}
Given $1 < q < \8.$ Since $M^{(q)}_{\g}(f)=(M_{\g}(|f|^q))^{1/q},$ $M^{(q)}_{\g}$ is bounded on $L^{p}_{\alpha}(\mathbb{B}_n)$ for each $q < p \le \infty.$
\end{remark}

Let $1\le p <\infty$ and let $E$ be a complex Banach space. We write
$L^p_{\alpha}(\mathbb{B}_n, E)$ for the Banach space of strongly measurable
$E$-valued functions on $\mathbb{B}_n$ such that
\be
\left(\int_{\mathbb{B}_n}\|f(z)\|^p_E dv_{\alpha}(z)\right)^{\frac{1}{p}}<
\infty.
\ee
Recall that $f : \mathbb{B}_n \mapsto E$ is said to be holomorphic if for each $x^* \in E^*,$ $x^* f$ is holomorphic in $\mathbb{B}_n.$ It is known that (for example, see \cite{Muj1986}) if $f$ is holomorphic in this weak sense, then it is holomorphic in the stronger sense that $f$ is the sum of a power series
\be
f (z) = \sum_{J \in \mathbb{N}^n_0} x_J z^J,\quad z \in \mathbb{B}_n,
\ee
where $x_J \in E.$ (As usual, $\mathbb{N}_0 = \mathbb{N} \cup \{0\}.$) The class of all such functions is denoted by $\mathcal{H}(\mathbb{B}_n, E).$ We write $\mathcal{A}^p_{\alpha}( \mathbb{B}_n, E)$ for the class of
weighted $E$-valued Bergman space of functions $f \in
\mathcal{H} (\mathbb{B}_n, E) \cap L^p_{\alpha} (\mathbb{B}_n, E).$

Then, by merely repeating the proof of the scalar case (e.g., Theorem 3.25 in \cite{Zhu2005}), we have the following interpolation result.

\begin{lemma}\label{le:EvaluedInterp}
Let $E$ be a complex Banach space. Suppose $\alpha>-1$ and $1 \le p_0 < p_1 < \8.$ If
\be
\frac{1}{p}=\frac{1-\theta}{p_0} + \frac{\theta}{p_1}
\ee
for some $0 < \theta < 1,$ then
\be
\left [ \mathcal{A}^{p_0}_{\alpha}(\mathbb{B}_n, E), \A^{p_1}_{\alpha} (\mathbb{B}_n, E) \right ]_{\theta} = \mathcal{A}^p_{\alpha}(\mathbb{B}_n, E)
\ee
with equivalent norms.
\end{lemma}

Now we are ready to prove Theorem \ref{th:tentspace}.

\

{\it Proof of Theorem \ref{th:tentspace}}.\;
By Lemmas \ref{le:BergmanballVolumeEsti}, \ref{le:BasicEsti01}, and \ref{le:BergmanFunctEsti01}, we have, for any $f \in \mathcal{H} (\mathbb{B}_n),$
\be
|f(z)|\lesssim \left ( \frac{1}{(1-|z|^2)^{n+1+\alpha}} \int_{D(z,\g)}|f|^q d v_{\alpha} \right)^{1/q} \lesssim A^{(q)}_{\g} (f)(z) \le M_{\g}^{(q)} (f)(z),
\ee
for all $z \in \mathbb{B}_n.$ Then we have that (3) implies (2), (2) implies (1) in Theorem \ref{th:tentspace}. By Remark \ref{rk:Maximalfunction} we conclude that (1) implies (3) when $q<p<\infty.$ Then, we need only to prove that for each $0 < p \le q,$
\be
\| M_{\g}^{(q)} (f) \|_{p, \alpha} \lesssim \|f \|_{p, \alpha}, \quad \forall f \in \A^p_{\alpha}.
\ee
On the other hand, by Lemmas \ref{le:BergmanballVolumeEsti} and \ref{le:BasicEsti01}, we know that
\be
M^{(q)}_{\g} (f)(z) \lesssim A^{(q)}_{2 \g} (f)(z),\quad \forall z \in \mathbb{B}_n.
\ee
It then remains to show that $\|A^{(q)}_{\g} (f)\|_{p,{\alpha}}\lesssim\|f\|_{p,{\alpha}}$ for $0 < p \le q.$

At first, we prove the case $0<p \le 1.$ For any $f \in \mathcal{A}^p_{\alpha}$ with the atomic decomposition \eqref{eq:atomdecomp}, one has by Lemma \ref{le:BasicEsti02} \be\begin{split}
A^{(q)}_{\g} (f) (z) = & \left ( \int_{D(z,\g)}|f|^q d \tau \right )^{\frac{1}{q}}\\
\le & \sum^{\infty}_{k=1}|c_{k}|\left(\int_{D(z,\g)}\frac{(1-|a_k|^2)^{q(pb-n-1-\alpha)/p}}{|1-\langle
w,a_{k}\rangle|^{ q b}}d\tau\right)^{\frac{1}{q}}\\
\le & \sum^{\infty}_{k=1}|c_{k}|\frac{(1-|a_k|^2)^{(pb-n-1-\alpha)/p}}{|1-\langle
z,a_{k}\rangle|^{b}}.
\end{split}\ee
Hence, by Lemma \ref{le:JfunctEsti},
\be\begin{split}
\int_{\mathbb{B}_n} |A^{(q)}_{\g} (f) |^p d v_{\alpha}
\le&\sum^{\infty}_{k=1}|c_{k}|^p(1-|a_k|^2)^{pb-n-1-\alpha}\int_{\mathbb{B}_n}\frac{1}{|1-\langle
z,a_{k}\rangle|^{pb}}dv_{\alpha}\\
\lesssim&\sum^{\infty}_{k=1}|c_{k}|^p.
\end{split}\ee
This concludes that
\be
\int_{\mathbb{B}_n} |A^{(q)}_{\g} (f) |^p d v_{\alpha} \lesssim \inf \Big \{ \sum_{k=1}^{\infty} |c_{k}|^{p} \Big \} \lesssim \|f\|^p_{L^p_{\alpha}}.
\ee

For the remaining case $1<p\le q$ we adopt the interpolation. Set $E = L^q(\mathbb{B}_n, \chi_{D(0, \gamma)} d \tau )$ and consider the operator
\be
T (f)(z,w) = f( \varphi_z (w)), \quad f \in \mathcal{H} ( \mathbb{B}_n).
\ee
Note that $\varphi_z (D(0, \gamma)) = D(z, \gamma)$ and the measure $d\tau$ is invariant under any automorphism of $\mathbb{B}_n$ (cf. Proposition 1.13 in \cite{Zhu2005}), we have
\be\begin{split}
\| T (f) (z) \|_E & = \left ( \int_{\mathbb{B}_n} \big | f  ( \varphi_z (w)) \big |^q \chi_{D(0, \gamma)} (w) d\tau (w) \right )^{\frac{1}{q}}\\
& = \left ( \int_{\mathbb{B}_n}|f ( w )|^q \chi_{D(z, \gamma)} (w) d\tau (w) \right )^{\frac{1}{q}}\\
& = A^{(q)}_{\g} (f)(z).
\end{split}\ee
On the other hand,
\be\begin{split}
\|T f\|^q_{L^q_{\alpha}(\mathbb{B}_n, E)}=&\int_{\mathbb{B}_n}\| T (f) (z) \|^q_E dv_{\alpha}\\
=&\int_{\mathbb{B}_n}|A^{(q)}_{\gamma}(f)(z)|^q dv_{\alpha}\\
\thickapprox & \int_{\mathbb{B}_n}|f|^q d v_{\alpha}.
\end{split}\ee
Then, we conclude that $T$ is bounded from $\mathcal{A}^q_{\alpha}$ into $L^q_{\alpha}(\mathbb{B}_n, E).$ Thus, combining with the case of $p=1$ proved above, we conclude from Lemma \ref{le:EvaluedInterp} that $T$ is bounded from $\mathcal{A}^p_{\alpha}$ into $L^p_{\alpha} (\mathbb{B}_n, E)$ for any $1<p<q,$ that is,
\be
\| A^{(q)}_{\g}(f) \|_{p, \alpha} \le C \| f \|_{p,\alpha},\quad \forall f \in \mathcal{A}^p_{\alpha},
\ee
where $C$ depends only on $\g, n, p, q$ and $\alpha.$ This completes the proof of Theorem \ref{th:tentspace}.\hfill $\Box$


\section{Proof of Theorem \ref{th:g-FunctCharat}}\label{proofg-funct}

In this section, we will prove Theorem \ref{th:g-FunctCharat}. Since
\be
(1- |z|^2) | \mathcal{R} f (z) | \le (1- |z|^2) |\nabla f (z) | \le |\tilde{\nabla} f (z)|, \quad \forall z \in \mathbb{B}_n,
\ee
we have that (4) implies (3), and (3) implies (2) in Theorem \ref{th:g-FunctCharat}. It remains to show that (2) implies (1) and (1) implies (4).

{\it Proof of $(2) \Rightarrow (1).$}\;
Recall that for any $f \in \mathcal{H}^p (\mathbb{S}_n)$ (the holomorphic Hardy space on $\mathbb{S}_n$),
\be
\int_{\mathbb{S}_n} | f (\zeta) - f(0) |^p d \sigma (\zeta) \thickapprox \int_{\mathbb{S}_n} \Big ( \int^1_0 |(1-s) \mathcal{R} f (s \zeta) |^q \frac{d s}{1 - s} \Big )^{\frac{p}{q}} d \sigma (\zeta).
\ee
(e.g., see \cite{AB1988}.) Then, noting that $f_{r}(z)=f(rz),$
\be\begin{split}
\int_{\mathbb{B}_n} & |f(z)-f(0)|^pdv_{\alpha}(z) \\
= & 2n \int_0^1 \int_{\mathbb{S}_n}|f(r\zeta)-f(0)|^p d\sigma(\zeta) r^{2n-1}(1-r^2)^{\alpha}dr\\
\thickapprox & 2n \int_0^1\int_{\mathbb{S}_n} \Big ( \int^1_0 |(1-s) \mathcal{R} f_r (s \zeta) |^q \frac{d s}{1 - s} \Big )^{\frac{p}{q}} d\sigma(\zeta) r^{2n-1}(1-r^2)^{\alpha} d r\\
= &  \int_{\mathbb{B}_n} \Big ( \int^1_0 |(1-s) \mathcal{R} f (s z ) |^q \frac{d s}{1 - s} \Big )^{\frac{p}{q}} d v_{\alpha} (z).
\end{split}\ee
This completes the proof that $(2) \Rightarrow (1).$

\

{\it Proof of $(1) \Rightarrow (4).$}\; We first consider the case $0<p\leq 1.$ To this end, we write
\be
f_k(z)= \frac{(1-|a_k|^2)^{(pb-n-1-\alpha)/p}}{(1-\langle
z,a_{k}\rangle)^b}.
\ee
An immediate computation yields that
\be
\nabla f_k(z)=\frac{ b \overline{a}_k(1-|a_k|^2)^{(pb-n-1-\alpha)/p}}{(1-\langle
z,a_{k}\rangle)^{b+1}}
\ee
and
\be
\mathcal{R} f_{k}(z)=\frac{b\langle z,a_k\rangle(1-|a_k|^2)^{(pb-n-1-\alpha)/p}}{(1-\langle
z,a_{k}\rangle)^{b+1}}
\ee
Then we have
\be\begin{split}
|\tilde{\nabla} f_k(rz)|^2=&(1-|rz|^2)(| \nabla f_k (rz) |^2 - | \mathcal{R} f_k (rz)|^2 )\\
=& b^2(1-|rz|^2)(1-|a_k|^2)^{2(pb-n-1-\alpha)/p}\frac{|a_k|^2-|\langle rz,a_{k} \rangle|^2}{|1-\langle
rz,a_{k}\rangle|^{2(b+1)}}
\end{split}\ee

Notice that
\beq\label{eq:BasicEsti03}
|1-t \lambda|\approx (1-t)+|1- \lambda|,\quad 0\leq t \leq 1, \; 0 \le |\lambda| \le 1.
\eeq
Set $\varepsilon_0 =\frac{(\alpha+1)q}{2p}.$ Then by \eqref{eq:BasicEsti03} one has
\be\begin{split}
G^{(q)}_{\tilde{\nabla} }(f_k) (z) = & b(1-|a_k|^2)^{(pb-n-1-\alpha)/p}\\
& \quad \times \left ( \int_0^1 \frac{(1-|rz|^2)^{q/2} (|a_k|^2-|\langle rz,a_{k} \rangle|^2)^{q/2}}{|1-\langle
r z,a_{k}\rangle|^{q (b+1)}}\frac{d r}{1-r |z|} \right )^{\frac{1}{q}}\\
\lesssim &  (1-|a_k|^2)^{(pb-n-1-\alpha)/p}\\
& \quad \times  \left( \int_0^1\frac{dr}{|1- r \langle z,a_{k}\rangle|^{q b+q/2}(1-r |z|)^{1-q/2}}\right )^{\frac{1}{q}}\\
\lesssim & (1-|a_k|^2)^{(pb-n-1-\alpha)/p} \\
& \quad \times \left( \int_0^1\frac{dr}{|1- r \langle
z,a_{k}\rangle|^{q b-\varepsilon_0} (1-r |z|)^{1+\varepsilon_0}}\right)^{\frac{1}{q}}\\
\lesssim & (1-|a_k|^2)^{(pb-n-1-\alpha)/p} \\
& \quad \times \left( \int_0^1\frac{dr}{\big[(1- r)+|1- \langle
z,a_{k}\rangle|\big]^{q b-\varepsilon_0}(1- r |z|)^{1+\varepsilon_0}}\right)^{\frac{1}{q}}\\
\lesssim & (1-|a_k|^2)^{(pb-n-1-\alpha)/p} \\
& \quad \times \left( \int_0^1\frac{dr}{|1- \langle
z,a_{k}\rangle|^{q b-\varepsilon_0}\big[(1- r)+(1-|z|)\big]^{1+\varepsilon_0}}\right)^{\frac{1}{q}}\\
\lesssim &  (1-|a_k|^2)^{(pb-n-1-\alpha)/p}\frac{(1-|z|)^{-\varepsilon_0 /q}}{|1-\langle
z,a_{k}\rangle|^{b-\varepsilon_0 /q}}.\\
\end{split}\ee
Hence, for any $f \in \mathcal{A}^p_{\alpha}$ with the atomic decomposition \eqref{eq:atomdecomp},
\be\begin{split}
\int_{\mathbb{B}_n}| & G^{(q)}_{\tilde{\nabla}}( f)(z)|^p d v_{\alpha}\\
& \leq \sum_{k=1}^{\infty} |c_k|^p \int_{\mathbb{B}_n}|G^{(q)}_{\tilde{\nabla}} (f_{k})(z)|^p d v_{\alpha}\\
& \lesssim  \sum_{k=1}^{\infty} |c_k|^p(1-|a_k|^2)^{pb-n-1-\alpha}\int_{\mathbb{B}_n}\frac{(1-|z|^2)^{\alpha-{\varepsilon_0 p}/q}}{|1-\langle z,a_{k}\rangle|^{pb-{\varepsilon_0 p}/q}}dv(z)\\
& \lesssim  \sum_{k=1}^{\infty}|c_k|^p(1-|a_k|^2)^{pb-n-1-\alpha}(1-|a_k|^2)^{-pb+n+1+\alpha}\\
& \lesssim  \sum_{k=1}^{\infty}|c_k|^p,
\end{split}\ee
the last second inequality is obtained by Lemma \ref{le:JfunctEsti} (notice that $\alpha-{\varepsilon_0 p}/q>-1$). The proof of the case $0<p\leq1$ is complete.

Now we turn to consider the case $1<p<\infty.$ In this case, $f$ has the integral representation
\be
f(z)=\int_{{\mathbb{B}}_n}\frac{f(w)dv_\alpha(w)}{(1-\langle z,w\rangle)^{n+1+\alpha}}
\ee
for all $z\in{\mathbb{B}}_n$. By Lemma \ref{le:BergmanFunctEsti02} one has
\be
\left|\tilde{\nabla}  f (r z) \right|\leq\sqrt{2}(n+1+\alpha)(1-|r z|^2)^{\frac{1}{2}}\int_{{\mathbb{B}}_n}\frac{|f(w)|(1-|w|^2)^{\alpha}dv(w)}{|1-\langle r z,w\rangle|^{n+1+\alpha+\frac{1}{2}}}.
\ee
Given $\varepsilon_0 =\frac{(\alpha+1)q}{2p}$ as above. Then one has by \eqref{eq:BasicEsti03} again
\be\begin{split}
G^{(q)}_{\tilde{\nabla}}( f) (z) =& \left(\int_0^1| \tilde{\nabla} f (r z) |^q \frac{d r}{1-|r z|}\right)^{1/q}\\
\lesssim&\left(\int_0^1(1-|rz|)^{q/2-1}\left|\int_{{\mathbb{B}}_n}\frac{|f(w)|(1-|w|^2)^{\alpha}dv(w)}{|1-\langle r z,w\rangle|^{n+1+\alpha+\frac{1}{2}}}\right|^qdr\right)^{1/q}\\
\le & \int_{{\mathbb{B}}_n}\left(\int_0^1\frac{(1-|rz|)^{q/2-1}|f(w)|^{q}(1-|w|^2)^{q\alpha}}{|1-\langle rz,w\rangle|^{q(n+1+\alpha+\frac{1}{2})}}dr\right)^{1/q}dv(w)\\
\lesssim & \int_{{\mathbb{B}}_n}|f(w)|(1-|w|^2)^{\alpha}\\
&\quad \times \left(\int_0^1\frac{(1-|rz|)^{q/2-1}}{\big[(1-r)+|1-\langle z,w\rangle|\big]^{q(n+1+\alpha)-\varepsilon_0}(1-r|z|)^{\varepsilon_0+\frac{q}{2}}}dr\right)^{1/q}dv(w)\\
\end{split}\ee
\be\begin{split}
\lesssim & \int_{{\mathbb{B}}_n}\frac{|f(w)|(1-|w|^2)^{\alpha}}{|1-\langle z,w\rangle|^{(n+1+\alpha)-\varepsilon_0 /q}}
\left ( \int_0^1\frac{dr}{\big[(1-r)+(1-|z|)\big]^{\varepsilon_0 + 1}}\right)^{1/q}dv(w)\\
\lesssim & \left ( \frac{1}{\varepsilon_0} \right )^{\frac{1}{q}}(1-|z|)^{-\varepsilon_0 /q}\int_{{\mathbb{B}}_n}\frac{|f(w)|(1-|w|^2)^{\alpha}}{|1-\langle z,w\rangle|^{(n+1+\alpha-\varepsilon_0 /q)}} d v(w).\\
\end{split}\ee
Thus, by Lemma \ref{le:BergmanOperEsti}, we have
\be
\| G^{(q)}_{\tilde{\nabla}}( f) \|_{p, \alpha} \lesssim \| f \|_{p, \alpha},
\ee
for $1<p<\infty.$ The proof of Theorem \ref{th:g-FunctCharat} is complete.

\section{Besov spaces}\label{Besov}

In this section we give an application of our main results (Theorems \ref{th:tentspace} and \ref{th:g-FunctCharat}) to obtaining new maximal and area integral characterizations of Besov spaces.

For $0 < p < \8$ and $-\8 < \alpha < \8$ we fix a nonnegative integer $k$ with $pk + \alpha > -1$ and define the so-called Bergman space $\mathcal{A}^p_{\alpha}$ introduced in \cite{ZZ2008} as the space of all $f \in \mathcal{H} (\mathbb{B}_n)$ such that $(1-|z|^2)^k \mathcal{R}^k f \in L^p (\mathbb{B}_n, d v_{\alpha}).$ One then easily observes that $\mathcal{A}^p_{\alpha}$ is independent of the choice of $k$ and consistent with the traditional definition when $\alpha > -1.$ Let $N$ be the smallest nonnegative integer
such that $pN + \alpha > -1$ and define
\beq\label{eq:BergmanSpaceNorm}
\| f \|_{p, \alpha} = | f(0)| + \left ( \int_{\mathbb{B}_n} (1-|z|^2)^{pN} | \mathcal{R}^N f (z) |^p d v_{\alpha} (z) \right )^{\frac{1}{p}},\quad f \in \mathcal{A}^p_{\alpha}\;.
\eeq
Equipped with \eqref{eq:BergmanSpaceNorm}, $\mathcal{A}^p_{\alpha}$ becomes a Banach space when $p \ge 1$ and a quasi-Banach space for $0 < p < 1.$

It is known that the family of the generalized Bergman spaces $\mathcal{A}^p_{\alpha}$ covers most of the spaces of holomorphic functions in the unit ball of $\mathbb{C}^n,$ such as the classical diagonal Besov space $B^s_p$ and the Sobolev space $W^p_{k,\beta}$ (e.g., \cite{BB1989}), which has been extensively studied before in the literature under different names (e.g., see \cite{ZZ2008} for an overview). We refer to Arcozzi-Rochberg-Sawyer \cite{ARS2006, ARS2008}, Tchoundja \cite{Tch2008} and Volberg-Wick \cite{VW} for some recent results on such Besov spaces and more references.

We have new maximal and area integral characterizations for these spaces as follows.

\begin{corollary}\label{cor:AreaMaximalCharat}
Suppose $\gamma >0, 1 \le q < \8,$ and $\alpha \in \mathbb{R}.$ Let $0 < p < \8$ and $k$ be a nonnegative integer such that $p k + \alpha > -1.$ Then for any $f \in \mathcal{H} (\mathbb{B}_n),$ $f \in \mathcal{A}^p_{\alpha}$ if and only if $A^{(q)}_{\gamma} (\mathcal{R}^k f)$ is in $L^p (\mathbb{B}_n, d v_{\alpha})$ if and only if $M^{(q)}_{\gamma} (\mathcal{R}^k f)$ is in $L^p (\mathbb{B}_n, d v_{\alpha}),$
where
\beq\label{eq:kareaFunct}
A^{(q)}_{\gamma} (\mathcal{R}^k f) (z) = \left ( \int_{D(z,\gamma)} \big | (1-|w|^2)^k \mathcal{R}^k f (w) \big |^q d \tau (w)  \right )^{\frac{1}{q}}
\eeq
and
\beq\label{eq:kmaximalFunct}
M^{(q)}_{\gamma} (\mathcal{R}^k f) (z)  = \sup_{z \in D(w, \g)}\left ( \int_{D(w,\gamma)} \big | (1-|u|^2)^k \mathcal{R}^k f (u) \big |^q  \frac{d v_{\alpha} (u)}{v_{\alpha} (D(w,\gamma))}\right )^{\frac{1}{q}}
\eeq

Moreover,
\beq\label{eq:kAreaFunctNorm}
\| f - f(0) \|_{p,\alpha} \approx \| A^{(q)}_{\gamma} (\mathcal{R}^k f) \|_{p, \alpha} \approx \| M^{(q)}_{\gamma} (\mathcal{R}^k f) \|_{p, \alpha},
\eeq
where ``$\approx$" depends only on $\gamma, q, \alpha, p, k,$ and $n.$
\end{corollary}

To prove Corollary \ref{cor:AreaMaximalCharat}, one merely notices that $f \in \mathcal{A}^p_{\alpha}$ if and only if $\mathcal{R}^k f \in L^p (\mathbb{B}_n, d v_{\alpha + pk})$ (e.g., Theorem 2.16 in \cite{Zhu2005}) and applies Theorem \ref{th:tentspace} to $\mathcal{R}^k f$ with the help of Lemma \ref{le:BasicEsti01}. When $\alpha > -1,$ we can take $k=0$ and then recover one of the main results in \cite{CO}.


Notice that $\mathcal{H}^p_s = \mathcal{A}^p_{\alpha}$ with $\alpha = -2 s -1,$ where $\mathcal{H}^p_s$ is the Hardy-Sobolev space defined as the set
\be
\left \{ f \in \mathcal{H} (\mathbb{B}_n ):\; \| f \|^p_{\mathcal{H}^p_s} = \sup_{0<r <1} \int_{\mathbb{S}_n} | (I + \mathcal{R} )^s f (r \zeta) |^p  d \sigma (\zeta) < \8 \right \}.
\ee
Here,
\be
(I + \mathcal{R} )^s f = \sum^{\8}_{k=0} (1 + k)^s f_k
\ee
if $f = \sum^{\8}_{k=0} f_k$ is the homogeneous expansion of $f.$ There are several real-variable characterizations of the Hardy-Sobolev spaces obtained by Ahern and Bruna \cite{AB1988}. These characterizations are in terms of maximal and area functions on the admissible approach region
\be
D_{\alpha} (\eta) = \left \{ z \in \mathbb{B}_n:\; | 1- \langle z, \eta \rangle | < \frac{\alpha}{2} (1 - |z|^2) \right \},\quad \eta \in \mathbb{S}_n,\; \alpha >1.
\ee
Now, Corollary \ref{cor:AreaMaximalCharat} presents new maximal and area integral descriptions of the Hardy-Sobolev spaces in terms of the Bergman metric. A special case of this is a characterization of the usual Hardy space $\mathcal{H}^p= \mathcal{A}^p_{-1}$ itself.

\subsection*{Acknowledgement} This research was supported in part by the NSFC under Grant No. 11171338.

\end{document}